\documentclass{amsart}
\usepackage{latexsym}
\usepackage{amssymb}
\usepackage{amsmath}
\usepackage{amsfonts}
\usepackage{graphics}
\newtheorem{theorem}{Theorem}
\newcommand{\R}{\ensuremath{\mathbb{R}}}

\newcommand{\bx}{\ensuremath{\mathbf{x}}}
\newcommand{\by}{\ensuremath{\mathbf{y}}}

\newcommand{\bp}{\ensuremath{\mathbf{p}}}

\newcommand{\br}{\ensuremath{\mathbf{r}}}

\newfont{\bb}{msbm10}

\def\RR{ {\mathbb R}}
\def\NN{ {\mathbb N}}

\def\sphere{{\mathbb S}}
\def\X{{\mathcal X}}
\def\P{{\mathcal P}}
\def\bfx{{\bf x}}
\def\bfy{{\bf y}}
\def\bff{{\bf f}}
\def\bfa{{\bf a}}
\def\half{{\frac{1}{2}}}

\def\bfalpha{{\mathbb \alpha}}
\def\bfbeta{{\mathbb \beta}}
\usepackage{graphicx}

\begin{document}

\title [Preconditioning for Approximation on the Sphere] 
{Stability and preconditioning for a hybrid approximation on the sphere}


\author{Q.~T.~Le~Gia}
\address{School of Mathematics and Statistics,
         University of New South Wales, Sydney, NSW 2052, Australia}
\email{qlegia@unsw.edu.au}	 
\author{Ian H.~Sloan}
\address{School of Mathematics and Statistics,
         University of New South Wales, Sydney, NSW 2052, Australia}
\email{i.sloan@unsw.edu.au}
\author{Andrew J.~Wathen}
\address{Mathematical Institute, Oxford University, 24-29 St Giles,
         Oxford OX1 3LB, UK}
\email{wathen@maths.ox.ac.uk}


\begin{abstract}
This paper proposes a new preconditioning scheme for a linear
system with a saddle-point structure arising from a hybrid approximation 
scheme on the sphere, an approximation scheme that combines (local) spherical
radial basis functions and (global) spherical polynomials. Making use
of a recently derived inf-sup condition \cite{SW08} and the Brezzi
stability and convergence theorem for this approximation scheme,
we show that the linear system can be optimally preconditioned with a
suitable block-diagonal preconditioner. Numerical experiments with a 
non-uniform distribution of data points support the theoretical conclusions.


\end{abstract}
\maketitle

\section{Introduction}\label{intro}

Amongst approaches for scattered data approximation on the sphere, 
the hybrid interpolation scheme of von Golitschek \& Light \cite{GL} 
and Sloan \& Sommariva \cite{SS06}, 
which employs both radial basis functions and spherical polynomials, 
seems an attractive 
method, especially when the data is concentrated in some regions (such as over mountain ranges and trenches on the Earth's surface), yet relatively sparse in other regions. 
The underlying idea is that radial basis functions can give good approximation for rapidly 
varying data over short distances, whereas the polynomial component can more effectively 
represent smooth variations on a global scale. The radial basis functions are centered 
at data points which are supposed given, and the linear combination of radial basis 
functions is constrained to be orthogonal, in a natural sense, to the finite dimensional 
space of polynomials.

However, the hybrid scheme poses difficulties in implementation,
compared with a pure radial basis function approximation, when the
number of centers is large. In the case of a pure radial basis function
approximation with a (strictly) positive definite kernel, the resulting
linear system has a matrix that is positive definite, allowing an 
iterative solution by the conjugate gradient method, and 
preconditioning by, for example, the additive Schwarz method,
see \cite{GST06}. For the hybrid scheme, in contrast, the 
linear equations 
for the relevant expansion
coefficients have the saddle-point form, see \cite{SS06},
\begin{equation}\label{saddle_general}
  \left[\begin{array}{cc}A &Q\\Q^T & 0\end{array}\right]
\left[\begin{array}{c}{\bfalpha}\\ {\bfbeta}\end{array}\right] =
\left[\begin{array}{c}{\bff}_X\\ {\bf{0}} \end{array}\right],\end{equation}
where $A\in\mathbb{R}^{N\times N}$ is a positive definite matrix arising from the radial basis function part of the function approximation, and $Q\in \mathbb{R}^{N\times M}$ is a matrix of 
spherical harmonics evaluated at the data centers, with $M\le N$. (The matrices $A$ and $Q$ are defined properly in Section 2).
The saddle-point structure means that the overall matrix is not positive
definite, and that the conjugate gradient method is no longer a suitable
iterative solver. More fundamentally, because the matrix has both
positive and negative eigenvalues, the problem of constructing a good 
preconditioner becomes more delicate. For a thorough review of strategies
and challenges for the numerical solution of saddle-point problems, 
see \cite{BenziGolubLiesen,ESW05}.

In this manuscript we concentrate on the stability of the saddle-point 
formulation
of this hybrid scheme, 
and devise and validate a rapid preconditioned iterative solution method for 
the solution of the equations from the approximation. We make use of the Brezzi stability and convergence theorem well known
in the context of mixed finite elements, along with the new inf-sup 
condition of \cite{SW08} to establish convergence of the approximation
scheme; and then use the inf-sup condition to obtain an optimal 
preconditioner.

A leading contender amongst the solution methods for equations with the 
structure \eqref{saddle_general} (see \cite{BenziGolubLiesen}) is the block
preconditioning method of \cite{MurGolWat00} employing an approximation 
to the Schur complement
\[
S:=Q^TA^{-1}Q.
\]
The use of
Schur complement approximations in preconditioners is by now well established (see for
example \cite {ESW05}) in the setting of mixed finite element methods. In particular, Verf\"urth \cite{Verf} showed that for the mixed finite element approximation of the Stokes flow problem the Schur complement is spectrally equivalent to the identity operator 
(or to the mass matrix or $L_2$ projection matrix in the finite element setting), 
making this a suitable approximation to the Schur complement. Verf\"urth's proof makes 
essential use of the Babuska-Brezzi or inf-sup condition (Assumption 2.1 in \cite{Verf}), 
see \cite{Brezzi74,Brezzi}.

In the present setting, the Schur complement turns out to be spectrally equivalent not to the identity operator/matrix, but rather to a specific diagonal but non-constant matrix.  This spectral equivalence, a main result of the paper, is stated in Theorem \ref{spectral_equiv}.
The key ingredient  here, in analogy with the known inf-sup condition for the Stokes flow problem,  is the inf-sup stability condition recently established by
Sloan and Wendland \cite{SW08} for the hybrid approximation problem.

An approximate solver for the primal operator
(the radial basis function interpolation matrix in this case) is
also required. This could be provided, for example,  by the domain decomposition
method of Le Gia, Sloan and Tran \cite{GST06}, or by any other preconditioner
for the pure radial basis function problem. The resulting block
diagonal preconditioner is symmetric and positive definite, hence
the preconditioned MINRES method (\cite{PS75},\cite{ESW05}) is
applicable to the full problem (which is symmetric but not positive
definite).

In Section 2 we formulate the hybrid approximation scheme and establish 
notation. Then in Section 3 we describe the inf-sup condition of \cite{SW08}, 
and use the Brezzi theorem to establish stability and convergence of the 
scheme. In Section 4 we turn to preconditioning, establishing there the 
main spectral 
equivalence result. In Section 5 numerical results are presented (using the
primal preconditioner of \cite{GST06}) and we conclude in Section 6.

\section{Problem formulation}\label{form}

Let $X=X_N=\{\bfx_1,\bfx_2,\ldots,\bfx_N\}$ be a set of $N$ distinct
points on the sphere $\sphere^d$ in $\RR^{d+1}$. Using these as centers we define
a radial basis function approximation space
$$\X_{X_N}=\X_N :=\{\sum_{i=1}^N \alpha_i \phi(\cdot,\bfx_i) : \alpha_1,\ldots,\alpha_N \in\RR\}$$
with a suitable kernel function $\phi$.
The kernel is assumed to be (strictly) positive definite, that is
$$\sum_{i=1}^N\sum_{j=1}^N\alpha_i\phi(\bx_i,\bx_j)\alpha_j\ge 0$$
for every set of points $X_N=\{\bx_1,\ldots,\bx_N\}\in \sphere^d$ and for all $N\in \NN$, with equality for distinct points $\bx_j$ only if $\alpha_1=\alpha_2= \ldots =\alpha_N=0.$

The native space $\mathcal{N}_{\phi}$ is defined as the completion under the inner product
\begin{equation}\label{IP}\langle \sum_{i} \alpha_i \phi(\cdot,\bfx_i),
\sum_{j} \alpha_j^{\prime} \phi(\cdot,\bfx_j)\rangle_{\phi} =
\sum_{i}\sum_{j}\alpha_i \alpha_j^{\prime}\phi(\bfx_i,\bfx_j)\end{equation}
of the linear space
$$F_{\phi} := \{\sum_{j=1}^N \alpha_j\phi(\cdot,\bfx_j),
\alpha_j\in\RR, \bfx_j\in\sphere^d, j=1,\ldots,N,N\in\NN\},$$
where we insist that the points $\bfx_j$ are distinct.
The norm is as usual defined by
$$\|\cdot\|_{\phi} = \langle \cdot,\cdot\rangle_{\phi}^{1/2} .$$
It is well-known that $\mathcal{N}_{\phi}$ is a reproducing kernel Hilbert space
(see \cite{Ar}) with the reproducing kernel $\phi(\cdot,\cdot)$. That is
\begin{eqnarray*}
\phi(\bfx,\bfy) &=& \phi(\bfy,\bfx),\quad \bfx,\bfy\in\sphere^d\\
\phi(\cdot,\bfy) &\in& \mathcal{N}_{\phi}, \quad \bfy\in\sphere^d\end{eqnarray*}
and for $f\in \mathcal{N}_{\phi}$
\begin{equation}\label{rep}\langle f,\phi(\cdot,\bfy)\rangle_{\phi} = f(\bfy), \quad \bfy\in\sphere^d .\end{equation}
The kernel function $\phi$ needs to be positive definite in
order that the inner product $\langle\cdot,\cdot\rangle_{\phi}$ satisfy the
positivity axiom for an inner product. Equivalently, 
the matrices $A_{X}$ defined by
\begin{equation}\label{AXphi}
\left( A_{X}\right)_{i,j} := \phi(\bfx_i,\bfx_j) ,\quad i,j = 1,\ldots,N
\end{equation}
are positive definite as well as symmetric for every $X$ and every $N\in\NN$.

Taking now a fixed $N\in\NN$ and a fixed set $X_N \subset \sphere^d$, we may define the usual radial basis function interpolant to a continuous function $f$ on $\sphere^d$ by
$$f_N(\bfx) = \sum_{j=1}^N \alpha_j \phi(\bfx,\bfx_j), $$
where $\alpha_1,\ldots,\alpha_N$ are such that
$$f_N(\bfx_i) = f(\bfx_i),\quad i= 1,\ldots,N ,$$
which we may write as
$$\sum_{j=1}^N\phi(\bx_i,\bx_j)\alpha_j=f(\bx_i) \textrm{ for }i= 1,\ldots,N. $$
That is, the vector
$\bfalpha = (\alpha_1,\ldots,\alpha_N)^T$ of coefficients satisfies
\begin{equation}\label{radb}
A_{X}\bfalpha = \bff_X,\end{equation}
where
\begin{equation}\label{fdef}\bff_X :=(f(\bfx_1),\ldots,f(\bfx_N))^T.\end{equation}

The hybrid approximation scheme of von Golitschek \& Light \cite{GL} and 
Sloan \& Sommariva, see \cite{SW08}, employs not only the radial basis 
functions, but also
spherical polynomials of total degree up to some conveniently chosen $L\ge 0.$  We define
$$\P_L = \hbox{span}\{Y_{\ell,k}:\, k=1,\ldots,M(d,\ell),\,\ell=0,\ldots,L\},$$
where $Y_{\ell,k}$ is a spherical harmonic of degree $\ell$, that is, the restriction to $\sphere^d$ of a homogeneous harmonic
polynomial in $\RR^{d+1}$ of degree $\ell$, and $M(d,\ell)$ is the dimension of the space spanned by the spherical harmonics of degree $\ell$.  Then $\P_L$ is the set of spherical polynomials of degree $\le L$.
We shall assume that $\{Y_{\ell,k}:\, k=1,\ldots,M(d,\ell),\,\ell=0,1,\ldots\}$ is an orthonormal set with respect to the usual $L_2$ inner product, that is
\[ \int_{\sphere^d}Y_{\ell,k}(\bfx)Y_{\ell',k'}(\bfx)d\omega(\bfx)=\delta_{\ell,\ell'}\delta_{k,k'},\]
where $d\omega(\bfx)$ denotes surface measure on $\sphere^d$.
Then it is well known that 
$\{Y_{\ell,k}:\, k=1,\ldots,M(d,\ell),\,\ell=0,1,\ldots\}$ is a complete orthonormal basis for
$L_2(\sphere^d)$.

For a given function $f$, the hybrid approximation scheme
is then to find
\begin{equation}\label{unl}
u_{N,L}(\bfx)= \sum_{j=1}^N \alpha_j \phi(\bfx,\bfx_j) \in\X_N
\end{equation}
and
\begin{equation}\label{pnl}
p_{N ,L}(\bfx)= \sum_{\ell=0}^L \sum_{k=1}^{M(d,\ell)} \beta_{l,k} Y_{\ell,k}(\bfx)\in\P_L
\end{equation}
such that
\begin{equation}\label{interp}
u_{N,L}(\bx_i)+p_{N,L}(\bx_i)=f(\bx_i),
\end{equation}
or equivalently,
\begin{equation}\label{1st}
\sum_{j=1}^N \alpha_j \phi(\bfx_i,\bfx_j) +
  \sum_{\ell=0}^L \sum_{k=1}^{M(d,\ell)} \beta_{\ell,k} Y_{\ell,k}(\bfx_i) = f(\bfx_i),
\quad i= 1,\ldots, N ,\end{equation}
which is to be solved subject to the side condition
\begin{equation}\label{2nd}
\sum_{j=1}^N \alpha_j q(\bfx_j) = 0\quad \forall q\in\P_L.\end{equation}
The condition (\ref{2nd}) is equivalent, via \eqref{rep}, to $\langle q, u_{N,L}\rangle_\phi=0$ for all $q\in\P_L$, forcing the radial basis function component to be $\mathcal{N}_\phi$-orthogonal to $\P_L$.  It also ensures that the defining linear system is square and symmetric.

The conditions (\ref{1st}),(\ref{2nd}) can also be seen to be those which 
derive from the solution of the constrained optimization problem
$$
\min_{u_{N,L}\in \X_N}\;\,\frac{1}{2}\| u_{N,L} - f\|_\phi^2\;\; \mbox{subject to} \;\;
\langle Y_{\ell,k},u_{N,L}\rangle_\phi=0$$
for all $\ell=0,\ldots,L , k=1,\ldots,{M(d,\ell)}$, the coefficients $\beta_{\ell,k}$ being the
Lagrange multipliers in the Lagrangian
\begin{eqnarray*}
{\mathcal L} &=&\frac{1}{2}\| u_{N,L} - f\|_\phi^2 +
\sum_{\ell=0}^L \sum_{k=1}^{M(d,\ell)} \beta_{\ell,k}\langle Y_{\ell,k},u_{N,L}\rangle_\phi\\
&=& \frac{1}{2}\| u_{N,L} - f\|_\phi^2 + \langle p_{N ,L},u_{N,L}\rangle_\phi .
\end{eqnarray*}
This is therefore another way of expressing the hybrid approximation problem.


By choosing the spherical harmonic functions as the basis for $\P_L$ in (\ref{2nd}),
we can write
(\ref{1st}),(\ref{2nd}) as a so-called `saddle-point' linear system of equations
\begin{equation}\label{saddle}
  \left[\begin{array}{cc}A_{X} &Q_{X,L}\\Q_{X,L}^T & 0\end{array}\right]
\left[\begin{array}{c}{\bfalpha}\\ {\bfbeta}\end{array}\right] =
\left[\begin{array}{c}{\bff}_X\\ {\bf{0}} \end{array}\right],\end{equation}
where $A_{X}$ is defined by (\ref{AXphi}),
 $\bfalpha$
is the vector of coefficients $\alpha_j$ as defined above,  
$\bfbeta$ is a vector containing the
coefficients $\beta_{\ell,k}$ for $k=1,\ldots,M(d,\ell),\,\ell=0,\ldots,L$, 
and $Q_{X,L}$ is the $N \times M$ matrix defined by
\begin{equation}\label{Qdef}(Q_{X,L})_{i,\ell k} := Y_{\ell,k}(\bfx_i) ,
\quad i = 1,\ldots,N, \; k=1,\ldots,M(d,\ell), \; \ell=0,\ldots,L,
\end{equation}
and 
$$M := \sum_{\ell=0} M(d,\ell) = \mbox{dim }(\P_L).$$

In the present application we need to prescribe more precisely the nature of the kernel $\phi(\bx,\by)$.  In the first place we shall assume that it is zonal, meaning that
$$ \phi(\bx,\by)=\Phi(\bx\cdot\by)$$
for some function $\Phi\in C[-1,1]$, where $\bx\cdot\by$ denotes the Euclidean inner product in $\R^{d+1}$. More precisely, we shall assume that $\phi(\bx,\by)$ has an expansion of the form
\begin{equation}\label{ass2} \phi(\bx,\by)=\sum_{\ell=0}^\infty\sum_{k=1}^{M(d,\ell)}a_\ell Y_{\ell,k}(\bx)Y_{\ell,k}(\by),\end{equation}
with $a_\ell>0$ for all $\ell\ge0$.
That the expansion is zonal follows from the addition theorem for spherical harmonics,
\[\sum_{k=1}^{M(d,\ell)}Y_{\ell,k}(\bfx)Y_{\ell,k}(\bfy)=\frac{M(d,\ell)}{\omega_d}P_\ell(d+1,\bfx\cdot\bfy),\]
where $P_\ell(d+1,z)$ is the Legendre polynomial of degree $\ell$ in dimension $d+1$ normalized to $P_\ell(d+1,1)=1$, and $\omega_d$ is the total surface measure of $\sphere^d$,
\[\omega_d=\int_{\sphere^d}d\omega(\bfx).\]
In this situation it is well known that the inner product in $\mathcal{N}_\phi$ can be written as
\begin{equation}\label{fiprod}\langle u,v\rangle_\phi = \sum_{\ell=0}^\infty\sum_{k=1}^{M(d,\ell)}\dfrac{\widehat{u}_{\ell,k}\widehat{v}_{\ell,k}}{a_\ell},\end{equation}
where
$$\widehat{u}_{\ell,k}= \int_{\sphere^d}u(\bx)Y_{\ell,k}(\bx)d\omega(\bx).$$
Indeed, as a special case of (\ref{fiprod}) we find
$$\langle f,\phi(\cdot,\by)\rangle_\phi = \sum_{\ell=0}^\infty\sum_{k=1}^{M(d,\ell)}\dfrac{\widehat{f}_{\ell,k}a_\ell Y_{\ell,k}}{a_\ell}=\sum_{\ell=0}^\infty\sum_{k=1}^{M(d,\ell)}\widehat{f}_{\ell,k} Y_{\ell,k}(\by)=f(\bfy),$$
thus verifying the reproducing kernel property (\ref{rep}).

If we further assume that for large $\ell$
\begin{equation}\label{ass3}
a_\ell\sim(\ell+1)^{-2s},
\end{equation}
then it follows from (\ref{fiprod}) and (\ref{ass3}) that the native space $\mathcal{N}_\phi$ is equivalent to the Sobolev space $H^s(\sphere^d)$ with inner product
\begin{equation}\label{Hsnorm}\langle u,v\rangle_{H^s}=\sum_{\ell=0}^\infty\sum_{k=1}^{M(d,\ell)}(\ell+1)^{2s}\widehat{u}_{\ell,k}\widehat{v}_{\ell,k}.\end{equation}

Technicalities aside, we remark that the essential difficulty in analysing the hybrid approximation is that $\X_N$ and $\P_L$ are both subsets of $\mathcal{N}_\phi$, and that in an appropriate sense both can approximate $\mathcal{N}_\phi$ as $N$ or $L$ tend to $\infty$.  We need the inf-sup condition, now to be introduced, to allow both subspaces to coexist comfortably within the one approximation.

\section{Inf-sup condition, and the Brezzi theorem}\label{sec_infsup}

Typically, uniqueness of the solution and optimal error estimates for saddle-point problems
follow from so-called inf-sup conditions together with appropriate coercivity of the primal operator.  For us an essential tool will be the following inf-sup theorem proved in \cite{SW08}.
In this theorem $h_X$, for a given point set $X=X_N\subset\sphere^d$, is the mesh norm, defined by
$$h_X:=\sup_{\bx\in\sphere^d}\inf_{\bx_j\in X}\rm{cos}^{-1}(\bx\cdot\bx_j).$$
In words, $h_X$ is the maximum geodesic distance from a point on $\sphere^d$ to the nearest point of $X$.
\begin{theorem}\label{SW_theorem} 
Let $\phi$ be a kernel satisfying (\ref{ass2}) and (\ref{ass3}) for 
some $s>d/2$. There exist constants $\gamma>0$ and $\tau>0$ depending only on $d$ and $s$ such that for all $L\ge 1$ and all $X_N = \{\bx_1,\ldots,\bx_N\}\subset \sphere^d$ satisfying $h_X\le\tau/L$ the following inequality holds:

$$\inf_{p\in \P_L\setminus\{0\}}\sup_{v\in\X_N\setminus\{0\}}\dfrac{\langle p,v\rangle_\phi}{\|v\|_\phi\|p\|_\phi}\ge\gamma.$$

\end{theorem}

To use this to prove stability we start with the following well-known theorem from Brezzi \cite{Brezzi}, which is at the heart of most analyses of mixed finite elements.

\begin{theorem}\label{brezzi} Let $H$ and $J$ be real Hilbert spaces, $a(\xi_1,\xi_2)$
be a continuous bilinear form on $H\times H$ and $b(\psi,\xi)$ be a continuous bilinear form
on $J\times H$. Let $\{H_N: N\in \NN\}$ and $\{J_L : L\in \NN\}$  be
sequences of subspaces of $H$ and $J$ respectively. Set
$$K=\{\eta \in H : b(\theta,\eta)=0 \;\;\forall \theta\in J \},\;\; K_{N,L}=\{\eta \in H_N : b(\theta,\eta)=0 \;\;\forall \theta\in J_L \} .$$
If
\begin{equation}\exists \gamma_0 >0\; \text{ such that }\quad a(\eta,\eta) \geq \gamma_0\|\eta\|_H^2
\quad \forall \eta\in K \cup K_{N,L},\label{coerc}\end{equation}
and
$$\exists \gamma_1 >0\; \text{ such that }\quad \sup_{\eta\in H \setminus\{0\}}
\frac{b(\theta,\eta)}{\|\eta\|_H} \geq \gamma_1 \|\theta\|_J\quad\ \forall \theta\in J
$$
\begin{equation}\text{ and }\sup_{\eta\in H_N \setminus\{0\}}
\frac{b(\theta,\eta)}{\|\eta\|_H} \geq \gamma_1 \|\theta\|_J\quad\ \forall \theta\in J_L,
\label{infsup}\end{equation}
then for every $\ell_1\in H^{\prime}$ and  $\ell_2\in J^{\prime}$ and every $N,L>0$
the discrete problem of finding $\xi_{N,L}\in H_N$ and $\psi_{N,L}\in J_L$ such that
\begin{eqnarray*}
a(\xi_{N,L},\eta)+b(\psi_{N,L},\eta) &=& \langle \ell_1,\eta\rangle\quad\forall \eta\in H_N\\
b(\theta,\xi_{N,L})&=&\langle\ell_2,\theta\rangle\quad\forall\theta\in J_L
\end{eqnarray*}
has a unique solution, and there exists a constant $C=C(\gamma_0,\gamma_1)>0$
such that
$$\|\xi-\xi_{N,L}\|_H + \|\psi - \psi_{N,L}\|_J \leq C
\left(\inf_{\widehat{\xi}_N\in H_N} \|\xi-\widehat{\xi}_N\|_H +
\inf_{\widehat{\psi}_L\in J_L} \|\psi-\widehat{\psi}_L\|_J \right)$$
where $\xi\in H$ and $ \psi\in J$ are defined by
\begin{eqnarray*}
a(\xi,\eta)+b(\psi,\eta) &=& \langle \ell_1,\xi\rangle\quad\forall \eta\in H,\\
b(\theta,\xi)&=&\langle\ell_2,\theta\rangle\quad\forall\theta\in J.
\end{eqnarray*}
\end{theorem}

To apply this theorem, we first observe that the hybrid approximation with its defining equations  (\ref{1st}) and (\ref{2nd}) can be written, using the reproducing kernel property (\ref{rep}), as the problem of finding  $u_{N,L}\in \X_N$ and $p_{N,L} \in \P_L$ such that
\begin{equation}\label{firsts}\langle u_{N,L}, \eta\rangle_\phi+\langle p_{N,L},\eta\rangle_\phi=\langle f,\eta\rangle_\phi \quad\forall\eta\in\X_N,\end{equation}
\begin{equation}\label{seconds} \langle q,u_{N,L}\rangle_\phi=0\qquad \forall q\in\P_L.\end{equation}
To use Theorem \ref{brezzi} we take
$H=\mathcal{N}_{\phi}$, $J=\P_L$, $H_N=\X_N$ and $J_L=\P_L$, with
the inner product on $\mathcal{N}_{\phi}$ being defined by (\ref{IP}), and the bilinear forms $a(\cdot,\cdot)$ and $b(\cdot,\cdot)$
both equal to the $\mathcal{N}_{\phi}$ inner product.
The coercivity condition
(\ref{coerc}) is trivially satisfied on the whole space $H=\mathcal{N}_\phi$ with $\gamma_0=1$ since
$$a(u,u)=\langle u,u\rangle_{\phi} = \|u\|_{\phi}^2 .$$
The existence of a constant $\gamma_1$ independent of $N$ and $L$
satisfying (\ref{infsup}) is ensured by Theorem \ref{SW_theorem}, provided $h_X\le\tau/L$.
 The first part of Theorem \ref{brezzi} then confirms that the solution of the system (\ref{firsts}) and (\ref{seconds}) exists and is unique provided $h_X\le \tau/L$.  The last part of that theorem defines the comparison quantities:  it defines $u_L\in \mathcal{N}_\phi$ and $p_L\in \P_L$ such that
$$\langle u_L,\eta\rangle_\phi+\langle p_L,\eta\rangle_\phi=\langle f,\eta\rangle_\phi \quad\forall\eta\in \mathcal{N}_\phi,$$
$$ \langle q,u_L\rangle_\phi=0\qquad \forall q\in\P_L.$$
The second equation says that $u_L$ is orthogonal to the space $\P_L$. In principle the orthogonality is  with respect to the $\mathcal{N}_\phi$ inner product, but because of the zonal property of the kernel it is easy to see that this is the same as the $L_2$ orthogonal projection.
Indeed, from \eqref{fiprod} we have
\[\langle q,u_L\rangle_\phi=0\,
\forall q\in\P_L\implies(\widehat u_L)_{\ell,k}=0\mbox{ for }\ell\in[0,L]
\implies\langle q,u_L\rangle_{L_2}=0\,\forall q\in\P_L.\]
The first of the latter set of equations then becomes, on specialising the choice of $\eta$ to $q\in\P_L$,
$$\langle p_L,q\rangle_\phi=\langle f,q\rangle_\phi \quad\forall q\in\P_L,$$
thus $p_L$ is the orthogonal projection of $f$ on the subspace $\P_L$ with respect to either the $L_2$ or the $\mathcal{N}_\phi$ inner products.  We write this orthogonal projection as $P_L f$.  Now we can write $p_L=P_L f$, and $u_L=f-p_L=f-P_L f$. Theorem \ref{brezzi} with 
$\xi = u_L$ and $\psi = p_L$
now gives the following convergence result, recovering a result obtained by a direct argument in \cite{SW08}.  Note that even though we have taken $J=\P_L$ in the theorem, the constants $\gamma_0$ and $\gamma_1$ do not depend on $L$, and hence neither does $C$.
\begin{theorem}
Let $\phi$ be a kernel satisfying (\ref{ass2}) and (\ref{ass3}) for some $s>d/2$. There exist constants $C>0$ and $\tau>0$ depending only on $d$ and $s$ such that for all $L\ge 1$ and all $X = X_N = \{\bx_1,\ldots,\bx_N\}\subset \sphere^d$ satisfying $h_X\le\tau/L$ the solutions of \eqref{unl},\eqref{pnl} and \eqref{interp} satisfy
\begin{equation}\label{stability}
\|f-u_{N,L}-p_{N,L}\|_\phi 
\leq C
\inf_{\widehat{\xi}_N\in \X_N}
\|(f-P_Lf)-\widehat{\xi}_N\|_\phi.
\end{equation}
\end{theorem}

Explicit error bounds in the $L_2$ norm for $f\in H^s$ and $f\in H^{2s}$ can then be obtained as in \cite{SW08}.

\section{Preconditioning}\label{linear}
Now we turn our attention to the linear algebra aspects of the hybrid approximation described in Section 2.

We have noted already that the hybrid approximation can be written as the linear system (\ref{saddle}),
 with $A_X, Q_{X,L}$ and $\bff$ defined by (\ref{AXphi}),(\ref{Qdef}), and (\ref{fdef}).

The solution of saddle-point linear systems such as (\ref{saddle_general})
has received much attention in recent years - see
\cite{BenziGolubLiesen} for an overview of possible approaches. In
particular, it was shown in \cite{MurGolWat00} that
a suitable preconditioner for the saddle point system
\begin{equation}\label{saddle_mat}
  \left[\begin{array}{cc}A &Q\\Q^T   & 0\end{array}\right]
\end{equation}
with positive definite $A$ is
\begin{equation}\label{mgw}
  \left[\begin{array}{cc}A &0\\0   & S\end{array}\right],
\end{equation}
where

$$S = Q^T A^{-1} Q$$
is the Schur complement.  This is because of the remarkable fact that the product of (\ref{saddle_mat}) by the inverse of (\ref{mgw}) is a diagonalisable matrix with just three distinct eigenvalues, namely $1, (1\pm\sqrt{5})/2$. Thus an appropriate Krylov subspace
iteration such as MINRES (see \cite{PS75}) will converge in just
three iterations. While this is generally not a practical preconditioner, an approximate preconditioner of the form
\begin{equation}
  \left[\begin{array}{cc}\widehat{A} &0\\0   & \widehat{S}\end{array}\right],
\end{equation}
where $\widehat{A}$ is a preconditioner for the
problem (\ref{radb}) involving only $A$, and $\widehat{S}$
is a suitable Schur complement approximation, will lead to rapid convergence.

In the present work we shall assume that an approximate preconditioner $\widehat{A}$ for $A=A_X$ is already available; one example would be the domain decomposition preconditioner from \cite{GST06} - this is the one we employ in the numerical results presented in the next section.  Our interest here is in finding an appropriate approximation to the Schur complement $S_X=Q^T_{X,L}A_X^{-1}Q_{X,L}$.  We shall see that this is handed to us by the inf-sup result in Theorem \ref{SW_theorem}.

That inf-sup condition can be stated as
$$\sup_{v\in \X_N \setminus\{0\}}
\frac{\langle p,v \rangle_{\phi}}{\|v\|_{\phi}} \geq
\gamma_1 \|p\|_\phi\quad\ \forall p\in \P_L,$$
provided $h_X\le\tau/L$.  With the help of the Cauchy-Schwarz inequality, this can be strengthened to a two-sided inequality,
\begin{equation}
\|p\|_\phi = 
\sup_{v\in\mathcal{N}_\phi\setminus\{0\}}
\frac{\langle p,v \rangle_{\phi}}{\|v\|_{\phi}}
\ge\sup_{v\in \X_N \setminus\{0\}}
\frac{\langle p,v \rangle_{\phi}}{\|v\|_{\phi}} \geq
\gamma_1 \|p\|_\phi
\quad\ \forall p\in \P_L,  \label{sup}
\end{equation}
provided $h_X\le\tau/L$.
To find an equivalent matrix expression, we write $p\in\P_L$ and $v\in\X_N$ as
$$p= \sum_{\ell=0}^{L} \sum_{k=1}^{M(d,\ell)}\beta_{\ell,k}Y_{\ell,k}, \qquad v=\sum_{i=1}^N\alpha_i\phi(\cdot,\bx_i).$$
With the help of the reproducing property (\ref{rep}), we find
\begin{align*}\langle p,v\rangle_\phi & =\sum_{\ell=0}^{L} \sum_{k=1}^{M(d,\ell)}\sum_{i=1}^N\beta_{\ell,k}\alpha_i \langle Y_{\ell,k},\phi(\cdot,\bx_i)\rangle_\phi\\
& =\sum_{\ell=0}^{L} \sum_{k=1}^{M(d,\ell)}\sum_{i=1}^N\beta_{\ell,k}\alpha_i Y_{\ell,k}(\bx_i)=\bfbeta^TQ^T_{X,L}\bfalpha,
\end{align*}
and
$$\|v\|_\phi=\langle v,v\rangle_\phi^{1/2}=
(\sum_{i=1}^N\sum_{j=1}^N\alpha_i\alpha_j\phi(\bx_i,\bx_j))^{1/2}
=(\bfalpha^TA_X\bfalpha)^{1/2}.$$
Also, with the aid of 
(\ref{fiprod})
we obtain
$$\|p\|_\phi=\left(\sum_{\ell=0}^L\sum_{k=1}^{M(d,\ell)}\dfrac{ \beta_{\ell,k}^2}{a_\ell}\right)^{\half}
 =\left(\bfbeta^T\Lambda_L \bfbeta\right)^\half,$$
where $\Lambda_L$ is the $M\times M$ diagonal matrix given by
\begin{equation}\label{Lambda}
(\Lambda_L)_{\ell k,\ell^\prime k^\prime}=\delta_{\ell\ell^\prime}\delta_{kk^\prime}/a_{\ell}.
\end{equation}
Thus in matrix terms (\ref{sup}) can be written as
\begin{equation}
 (\bfbeta^T \Lambda_L \bfbeta)^{\half} \geq
\sup_{\bfalpha\in\RR^N\setminus\{0\}}
\frac{\bfbeta^T Q^T_{X,L}\bfalpha}{(\bfalpha^T A_{X}\bfalpha)^{\half}}
\geq \gamma_1 (\bfbeta^T \Lambda_L \bfbeta)^{\half}\quad \forall \bfbeta\in\RR^M.
\end{equation}
Because $A_{X}$ is symmetric and positive definite, the central term can be simplified by the substitution $\bfalpha=A_X^{-\half}\bfa$, making it
$$
\sup_{\bfa\in\RR^N\setminus\{0\}}
\frac{\bfbeta^T Q_{X,L}^T A_{X}^{-\half}\bfa}{(\bfa^T \bfa)^{\half}}\,=\,(\bfbeta^T Q_{X,L}^TA_X^{-1}Q_{X,L}\bfbeta)^\half,
$$
with the last step following because the supremum over $\bfa$ is clearly achieved by
$\bfa=(\bfbeta^T Q^T_{X,L} A_X^{-\half})^T$.
Thus in matrix terms (\ref{sup}) can be expressed as
\begin{equation}\label{schur}
\bfbeta^T \Lambda_L \bfbeta \geq \bfbeta^T Q^T_{X,L} A^{-1}_XQ_{X,L}\bfbeta\geq \gamma_1^2\bfbeta^T\Lambda_L\bfbeta\quad\forall\bfbeta\in\RR^M.
\end{equation}
Through the above arguments we have established the following theorem.

\begin{theorem}\label{spectral_equiv}
 Let $\phi$ be a kernel satisfying (\ref{ass2}) and (\ref{ass3}) for some $s>d/2$, and let $A_X$ and $Q_{X,L}$ be given by (\ref{AXphi}) and (\ref{Qdef}).  For all $L\ge 1$ and $h_X\le\tau/L$, where $\tau$ is as in Theorem \ref{SW_theorem}, the Schur complement $S_X=Q_{X,L}^T A_X^{-1}Q_{X,L}$ is spectrally equivalent to the diagonal matrix $\Lambda_L$ given by (\ref{Lambda}).
\end{theorem}

It follows from the theorem that our practical recommendation for the hybrid problem is a preconditioner of the form
\begin{equation}
  \left[\begin{array}{cc}\widehat{A} &0\\0   & \Lambda_L\end{array}\right],
\label{prec}\end{equation}
where $\widehat{A}$ is an approximation to $A$, and $\Lambda_L$ is defined by (\ref{Lambda}).
\section{Numerical examples}\label{numer}

\begin{figure}[ht]
\begin{center}
\scalebox{0.65}{\includegraphics{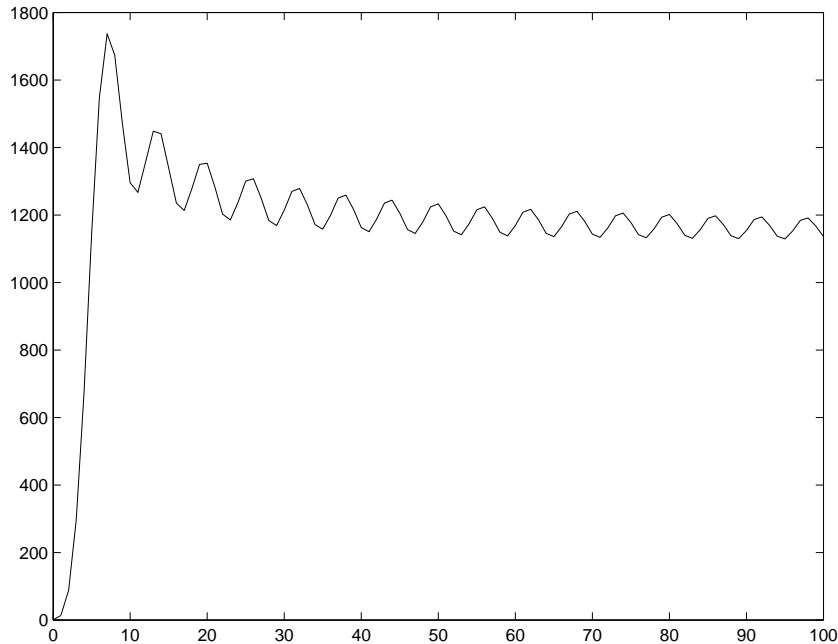}}
\caption{Numerical values of $(\ell+1)^5a_\ell$}\label{fig:aell}
\end{center}
\end{figure}

We will use the following kernel
\[
  \phi(\bx,\by) = \psi(|\bx-\by|) = \psi(\sqrt{2-2\bx \cdot \by}),
  \quad \bx,\by \in \sphere^2,
\]
where the radial basis function $\psi(r)$ is one of the three choices
\[
 \psi_0(r) = (1-r)^2_{+} , \quad \psi_1(r) = (1-r)^4_{+} (4r+1) , \quad 
\psi_2(r) = (1-r)^6_{+} (35 r^2+18r+3)
\]
with $(x)_{+} = x$ for $x\ge 0 $ and $0$ otherwise. Note that 
$\psi_0\in C^0(\RR^3), \psi_1\in C^2(\RR^3)$ and $\psi_2\in C^4(\RR^3)$ and each 
is positive definite (see \cite{Wend}).
We comment that we have not here employed any scaling of the compactly supported
RBFs. 
Using functions with smaller support would improve matrix
conditioning, but because it would also reduce approximation accuracy we have chosen 
not to use any scaling.
It is shown theoretically in \cite{NarWar02} and verified
numerically in Figure \ref{fig:aell}  
that the coefficients $a_\ell$ in the
expansion of the kernel $\phi$ defined from $\psi_1$ are of order $(1+\ell)^{-5}$. This
is consistent with (\ref{ass3}) but because the constants in the equivalence are large we 
have chosen to work directly with the Fourier-Legendre coefficients $a_{\ell}$,
which being $1$-dimensional integrals are easily evaluated numerically.
For the sphere $\sphere^2$, we have
\[
a_\ell = 2\pi \int_{-1}^1 \Phi(t) P_\ell(3,t) dt,
\quad \mbox{ where } \Phi(t) = \psi(\sqrt{2-2t}).
\]

For the preconditioner, as described in \cite{GST06}, the matrix
$A$ is preconditioned using a domain decomposition technique.
First, given a fixed parameter $0<\nu<\pi$, an appropriate set
of centers $\{\bp_1,\ldots,\bp_J\} \subset X$ is chosen over the
whole sphere so that
\[
   \min_{i\ne j} \cos^{-1} (\bp_i \cdot \bp_j) \ge \nu.
\]
Second, with another fixed parameter $0 < \mu< \pi/3$, we
decompose the point set $X$ into a collection of smaller sets
$X_j$, for $j=1,\ldots,J$, defined by
\[
  X_j := \{ \bx \in X : \cos^{-1}(\bx \cdot \bp_j) \le \mu \}.
\]
The sets $X_j$ with cardinality $m_j$, for $j=1,\ldots,J$, may overlap 
each other and must satisfy $\cup_{j=1}^J X_j = X$. The restriction operator from $\RR^N$ 
to $\RR^{m_j}$ is denoted by $R_j$ and the extension operator from $\RR^{m_j}$ 
to $\RR^N$ is $R^T_j$. Given a vector $\br \in \RR^N$, the
action of the preconditioner $\widehat{A}$ is given by 
\[
\widehat{A}^{-1} \br = \sum_{j=1}^J R^T_j (A_j)^{-1} R_j \br
\]
where the matrix $A_j$ is the restriction of the full matrix
$A$ on the subdomain $X_j$ (see \cite{GST06} for more details).

By Theorem \ref{spectral_equiv} the diagonal matrix $\Lambda_L$ defined by (\ref{Lambda}) 
is spectrally equivalent to $S_X = Q^{T}_{X,L} A^{-1}_X Q_{X,L}$.
The block diagonal preconditioner for (\ref{saddle_mat}) is therefore the matrix (\ref{prec}).

The results here are for interpolation of the function
\[
f(x,y,z) =  \exp(x+y+z)+ [0.01-x^2-y^2-(z-1)^2]^2_{+},
\]
consisting of a smooth first term and a second term whose support is a cap
of Euclidean radius $0.1$.

Using each of the kernel functions $\phi$ obtained from the radial basis functions 
$\psi_m, m=0,1,2$, we employ $N = 2000$, $4000$, $8000$, $16000$ and $32000$
points, and maximum polynomial degree $L = 0,5,10,15,20,25$. In each case, a thousand of the 
points were generated in a cap about the $z$ axis subtending an angle of $0.1$ radians at the 
origin and the remaining points distributed outside this cap. The Saff-Kuijlars equal area algorithm described in \cite{SaffK} was used to generate these points in the following manner. 
Firstly the Saff-Kuijlars points are generated on the whole sphere and those in the cap region are discarded. Then a similar equal area construction only for the cap is used to generate $1000$ points in this region. 

The number of MINRES iterations and the CPU time in seconds required for convergence 
to a residual norm tolerance of $10^{-9}$ are tabulated for the unpreconditioned case
in Table~\ref{tab:unprec} and for the preconditioner introduced here in Table~\ref{tab:prec}. The computer code is written in Fortran 90, compiled with the Intel compiler and run 
on a single core of an SGI Altix XE320 with two Intel Xeon X5472 CPUs.
\begin{table}
\begin{center}
\begin{tabular}{|c|l|r|r|r|r|r|}
\hline
$m$ & $N$    & 2000   &      4000  &     8000  &     16000  &   32000 \\
\hline
0 & $L= 0$ &  207 ( 14)&  228 ( 62)&  273 ( 297)& 294 ( 1229) &  367 ( 6277)\\
  & $L= 5$ &  813 ( 57)& 1116 (304)& 1531 (1662)&1996 ( 8297) & 2535 (43513)\\
  & $L=10$ & 1053 ( 77)& 1488 (413)& 2240 (2463)&3016 (12573) & 4001 (67005)\\
  & $L=15$ & 1080 ( 85)& 1492 (498)& 2208 (2475)&3037 (12788) & 4270 (71728)\\
  & $L=20$ & 1126 ( 99)& 1466 (464)& 2032 (2334)&2777 (17227) & 3832 (66025)\\
  & $L=25$ & 1262 (123)& 1467 (499)& 1976 (2382)&2562 (11069) & 3461 (58451)\\
\hline
1 & $L= 0$ & 1041 ( 75)&  790 ( 218)&  695 ( 755)&  912 ( 3954)&  995 (17436)\\ 
  & $L= 5$ & 2654 (193)& 4293 (1194)& 3564 (3889)& 3172 (13736)& 2321 (40754)\\
  & $L=10$ & 3647 (280)& 4370 (1237)& 4047 (4488)& 5335 (42445)& 3381 (59228)\\
  & $L=15$ & 3267 (265)& 4308 (1281)& 4021 (4594)& 3792 (16908)& 2808 (48999)\\
  & $L=20$ & 2883 (264)& 2923 ( 935)& 3329 (6906)& 3977 (17987)& $>$24 hours \\
  & $L=25$ & 2659 (265)& 2837 (1006)& 3295 (4048)& 2756 (12735)& 3461 (61280)\\
\hline
2& $L= 0$ & 1079 ( 82) & 1578 (457) &  1903 ( 2164) & 2743 (12233) & 2671 (48177)\\
 & $L= 5$ & 1978 (152) & 2594 (757) & 14424 (16476) & $>$24 hours & $>$40 hours \\  
 & $L=10$ & 2205 (176) & 3130 (929) & 13569 (15603) & $>$24 hours & $>$40 hours \\  
 & $L=15$ & 2402 (205) & 2803 (873) & 11970 (14053) & $>$24 hours & $>$40 hours \\  
 & $L=20$ & 1796 (169) & 1869 (624) &  8710 (10687) & 13551 (63016) &$>$40 hours\\
 & $L=25$ & 1615 (168) & 1758 (629) &  6417 ( 8239) & 10299 (48420) &$>$40 hours\\
\hline
\end{tabular}
\caption{MINRES iteration count (CPU time) without preconditioning}
\label{tab:unprec}
\end{center}
\end{table}

\begin{table}
\begin{center}
\begin{tabular}{|c|l|r|r|r|r|r|}
\hline
$m$ & $N$    & 2000   &   4000  &    8000 &     16000 &   32000 \\
\hline
0 &  $L=0$  & 31 ( 6)& 39 (18) & 30 ( 51)& 29 (225) &  39 ( 959)\\ 
  &  $L=5$  & 59 (11)& 71 (31) & 58 ( 96)& 62 (465) &  75 (1802)\\
  &  $L=10$ & 70 (13)& 88 (39) & 70 (116)& 71 (532) &  89 (2120)\\
  &  $L=15$ & 76 (15)& 93 (43) & 76 (128)& 76 (572) &  96 (2279)\\
  &  $L=20$ & 83 (17)& 98 (48) & 80 (138)& 82 (780) &  99 (2393)\\
  &  $L=25$ & 95 (20)& 97 (50) & 80 (142)& 84 (649) & 104 (2505)\\
\hline 
1 &  $L=0$  &  43 ( 8)&  75 (34)&  35 ( 66)&  29 (227)&  47 (1161) \\
  &  $L=5$  &  76 (15)& 128 (57)&  83 (138)&  74 (557)& 105 (2559)\\
  &  $L=10$ &  91 (18)& 148 (67)&  98 (163)&  94 (705)& 136 (3296)\\
  &  $L=15$ &  98 (19)& 168 (78)&  96 (163)& 100 (758)& 148 (3615)\\
  &  $L=20$ & 107 (22)& 170 (83)& 103 (180)& 103 (788)& 153 (3754)\\
  &  $L=25$ & 106 (23)& 174 (89)& 103 (185)& 113 (870)& 161 (3972)\\
\hline
2& $L= 0$ & 64 (13)  & 149 ( 61) & 46  ( 81) &  30 ( 243) & 61  (1419) \\  
 & $L= 5$ & 95 (19)  & 157 ( 70) & 88  (151) & 103 ( 797) & 140 (3380) \\  
 & $L=10$ & 95 (19)  & 171 ( 71) & 102 (176) & 111 ( 861) & 146 (3425) \\  
 & $L=15$ & 112 (23) & 187 (123) & 113 (199) & 118 ( 923) & 165 (3614) \\  
 & $L=20$ & 119 (25) & 197 ( 86) & 115 (207) & 133 (1052) & 196 (4301) \\  
 & $L=25$ & 125 (28) & 201 ( 91) & 119 (220) & 131 (1046) & 203 (4433) \\
\hline
\end{tabular}
\caption{MINRES iteration count(CPU time) with preconditioning
         $\widehat{A}$ and $\widehat{S}$}
\label{tab:prec}
\end{center}
\end{table}
The preconditioning is seen to be effective: as anticipated from the theory above, 
the number of iterations remains approximately constant over all choices of 
$N$ for each degree $L$. Indeed, for each $N$ aside from the simple case $L=0$ 
(in which the radial basis function approximation matrix is only 
supplemented by one row and one column),
the iteration counts grow only slowly 
with increasing $L$.

In order to determine the descriptiveness of the bound \eqref{schur}
we have also computed the generalised eigenvalues $\lambda_i$ of the pencil
$Q_{X,L}^T A_X^{-1} Q_{X,L} - \lambda\Lambda_L$ for the case $N=4000$ points, for $m=0,1$ 
and for $L=5,10,15,20,25$. Note that the generalised eigenvalues are exactly
the eigenvalues of $\Lambda^{-1}_L (Q_{X,L}^T A_X^{-1} Q_{X,L})$. 
The minimum and maximum computed eigenvalues
are given in Table~\ref{tab:eigvals}. It is noticeable how close the largest
eigenvalue is to the analytical upper bound of $1$ and that, although the lowest
eigenvalue does decrease for larger $L$, it remains reasonably
close to $1$. (Note that for fixed $X$ and increasing $L$ the inf-sup condition
must eventually break down.)
\begin{table}
\begin{center}
\begin{tabular}{|c|c|c|c|c|c|c|}
\hline
$m$ & $L$         & 5   & 10  & 15  & 20 &  25 \\
\hline
0 &$\lambda_{\min}$&0.9987434&0.9899326&0.9623012&0.9068357&0.8348191\\ 
  &$\lambda_{\max}$&0.9997653&0.9997658&0.9997674&0.9997753&0.9998099\\
\hline
1 &$\lambda_{\min}$&0.9999955&0.9999125&0.9993989&0.9973949&0.9908182\\ 
  &$\lambda_{\max}$&0.9999986&0.9999986&0.9999986&0.9999986&0.9999989\\
\hline  
\end{tabular}
\caption{Extreme eigenvalues of $Q_{X,L}^T A_X^{-1} Q_{X,L} - \lambda\Lambda_L$ for $N=4000$} \label{tab:eigvals}
\end{center}
\end{table}
\section{Conclusions}\label{concl}

By employing a recent inf-sup stability result of Sloan and Wendland,
we have derived an effective preconditioned iterative solver for the
hybrid radial basis function and spherical polynomial approximation
scheme of Sloan and Sommariva. The preconditioner requires only a
good approximation for the radial basis function interpolation
problem and a simple diagonal scaling matrix for the Schur complement 
based on the Fourier-Legendre coefficients of the kernel function used 
in the radial basis function interpolant. We have established theoretically 
that the preconditioner for the Schur complement is optimal. 

\vspace*{0.3in}

{\bf Acknowledgements.}
The support of the Australian Research Council 
is gratefully acknowledged.

\bibliographystyle{amsplain}

\end{document}